\documentclass[letterpaper, 11pt]{article}
\usepackage{amsmath, amsthm, amssymb}    
\usepackage{bridges}			
\usepackage{graphicx}			
\usepackage[colorlinks=true, urlcolor=blue, citecolor=black, linkcolor=black]{hyperref}  
\usepackage{subcaption}			
\usepackage{wrapfig}			
\usepackage{comment}

\newcommand{\R}{\mathbb R}
\newcommand{\Z}{\mathbb Z}
\renewcommand{\S}{\mathbb S}
\renewcommand{\H}{\mathbb H}
\newcommand{\Sol}{{\rm Sol}}

\usepackage{wrapfig}

\title{Non-Euclidean Virtual Reality IV: Sol}

\author{Remi Coulon\textsuperscript{1}, Elisabetta A. Matsumoto\textsuperscript{2}, Henry Segerman\textsuperscript{3}, Steve Trettel   \textsuperscript{4}
\vspace{10pt}\\
\begin{minipage}[b]{0.5\textwidth}
\centering
\textsuperscript{1}Univ Rennes, CNRS, France; remi.coulon@univ-rennes1.fr\\
\textsuperscript{3}Oklahoma State University; henry@segerman.org
\end{minipage}
\begin{minipage}[b]{0.5\textwidth}
\centering
\textsuperscript{2}Georgia Institute of Technology; sabetta@gatech.edu\\
\textsuperscript{4}Stanford University;\\ trettel@stanford.edu
\end{minipage}
}
\date{}					

\begin{document}

\maketitle

\thispagestyle{empty}

\begin{abstract}

This article presents virtual reality software designed to explore the Sol geometry.
The simulation is available on \href{http://www.3-dimensional.space/sol.html}{3-dimensional.space/sol.html}

\end{abstract}

\begin{figure}[h!tbp]
	\centering
	\includegraphics[width=0.8\textwidth]{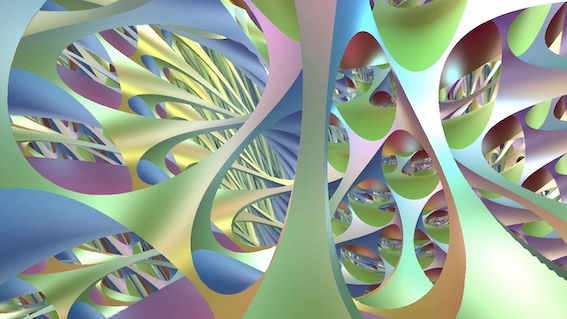}
		\caption{Intrinsic view of a Sol manifold built as a torus bundle over the circle with Anosov monodromy.}
	\label{fig:1}
\end{figure}

Geometrization (conjectured by Thurston, proved by Perleman) states that every closed three-dimensional manifold can be decomposed into elementary ``building blocks'' each of which is modeled onto a specific geometry.
There are eight such models: $\R^3$, $\S^3$, $\H^3$, $\H^2\times\R$, $\S^2\times\R$, Nil, Sol, and the universal cover of ${\rm SL}(2,\R)$.
We developed virtual reality software whose aim is to simulate these eight geometries.
We populate each of these metric spaces $X$ with various objects (spheres, planes, cylinders, lights, lattices, etc.) and compute what an observer would see if light follows the geodesics of $X$.
Using a virtual reality headset, the user can walk in these spaces and experience their surprising properties.
This paper presents an expository account of Sol geometry, and stems from a larger project to develop accurate, real time, intrinsic, and mathematically useful illustrations of homogeneous (pseudo)-riemannian spaces.

\section*{Sol Geometry}

\paragraph{Definition and metric.}
There are various ways to describe Sol geometry.
It can be seen as a Lie group $X$ whose underlying space is $\R^3$ with the group law given by 
\begin{equation*}
	(x_1,y_1,z_1) \ast (x_2,y_2,z_2) = (e^{z_1}x_2 + x_1, e^{-z_1}y_2 + y_1, z_1 + z_2)
\end{equation*}
The identity element is the point $o = (0,0,0)$ which we choose as the origin of the space.
Algebraically it is an extension of $\R^2$ by $\R$, and is therefore a \emph{solvable} group, hence the name of the geometry.
The space $X$ is endowed with a riemannian metric.
The metric tensor at an arbitrary point $p = (x,y,z)$ is given by 
\begin{equation*}
	ds^2 = e^{-2z}dx^2 + e^{2z}dy^2 + dz^2 \tag*{\cite[Section~1.7]{Troyanov:1998aa}}
\end{equation*}
With this metric, the action of $X$ on itself is an action by isometries.
In contrast to other geometries for which the underlying space is a group, $X$ has only finitely many more symmetries. 
These corresponds to the stabilizer of the origin; it acts by the rotations/reflections preserving the union of the $x,y$ axes.

\paragraph{Geodesic flow.}
In order to compute the trajectory of the light rays in Sol, we need a parametrization of its geodesics.
Using standard tools of riemannian geometry, one can prove that any curve $c \colon [0,1] \to X$ is a geodesic if and only if $c(t) = (x(t), y(t), z(t))$ satisfies the following differential equations
\begin{equation*}
	\left\{
	\begin{split}
		\ddot x & = 2\dot x\dot z  \\
		\ddot y & = - 2\dot y\dot z \\
		\ddot z & = - e^{-2z} \dot x^2 + e^{-2z} \dot y^2
	\end{split}
	\right.
\end{equation*}
This system can be explicitly solved using Jacobi's elliptic functions and Jacobi's zeta function \cite{Troyanov:1998aa}.
Nevertheless the explicit formulas are not particularly enlightening.
Instead, in the next sections we explore a few interesting features of the Sol geometry.

\begin{figure}[!htbp]
\centering
\includegraphics[width=0.7\textwidth]{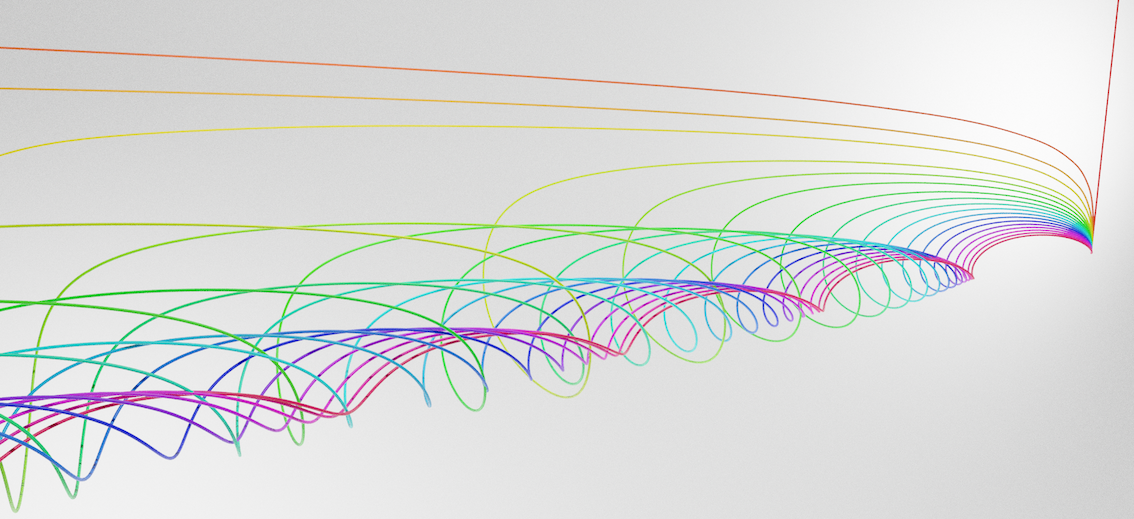}
\caption{Geodesics (for the Sol metric) starting at the origin in the  $\R^3$ model.}
\label{fig:2}
\end{figure}

\begin{figure}[!htbp]
	\centering
	\includegraphics[width=0.8\textwidth]{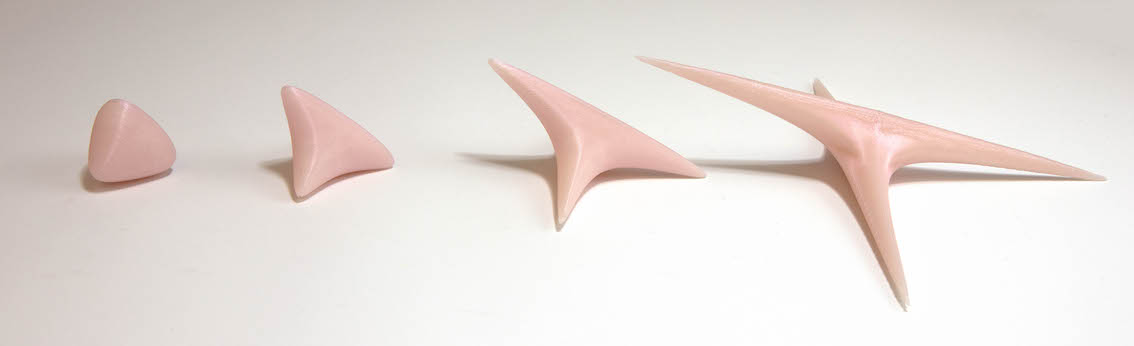}
	\caption{3D printed models of geodesic spheres in Sol.
	These models were a helpful tool to interpret the simulations described below. Photo: \copyright\ Edmund Harris}
	\label{fig:spheres}
\end{figure}

\paragraph{Shape of the Spheres.}
Spheres in Sol are quite surprising \cite{Coiculescu:2019aa}.
To interpret their shapes one needs to keep the following observation in mind.
Let $p = (x,y,z)$ be a point in $X$, with $z > 0$ and $c \colon [0,1] \to \R$ a geodesic from the origin $o$ to $p$.
Assume that one wants to extend the path $c$ to move further away from $o$.
Increasing the $x$-coordinate by $\epsilon$ will move us by a distance of $e^{-z}\epsilon$ (in the Sol metric).
On the other hand, increasing by $\epsilon$ the $y$-coordinate will move us by a much larger distance, namely $e^z\epsilon$.
Hence the upper part of the sphere (for $z > 0$) has a tendency to stretch along the $y$-axis, see \autoref{fig:spheres}.
The spheres exhibit the $D_8$-symmetry (where $D_8$ is the dihedral group of order $8$ and the stabilizer of the origin).
Hence the lower part of the sphere is stretched along the $x$-axis.

\section*{Exploration through Virtual Reality}

So far, we have reviewed the geometry of Sol \emph{extrinsically}, e.g. by referencing a model $(\R^3,ds^2)$.
We now explain what an observer living in a ``Sol world'' would see (the \emph{intrinsic} point of view).
To see in $\Sol$, we imagine our retina (or camera) as a small screen in space, and the image we see on the screen is given by tracing out along geodesics into the world, until they hit an object.

%
%
%
%
%
%
%
\paragraph{Hyperbolic foliation and horizontal planes.}
The Sol geometry has several remarkable two-dimensional subspaces.
Given $c \in \R$, the planes $\{x = c\}$ and $\{y = c\}$ are totally geodesic subspaces of Sol, which are isometric to the standard hyperbolic plane $\H^2$.
On the other hand, each plane $\{z = c\}$ is endowed with a distorted euclidean metric whose unit circle is an ellipse.
However it is not totally geodesic.
More precisely, the only geodesics contained in the plane $\{ z = 0\}$ are $\gamma_+(t) = (t, t, 0)$ and $\gamma_-(t) = (t, -t, 0)$.
This leads to some puzzling features.

\begin{figure}[h!tbp]
\centering
\begin{minipage}[b]{0.32\textwidth} 
	\includegraphics[width=\textwidth]{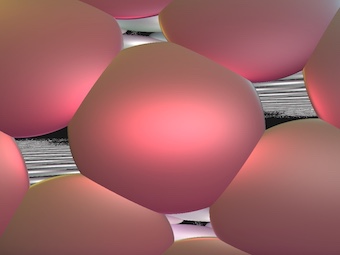}

\end{minipage}
~ 
\begin{minipage}[b]{0.32\textwidth} 
	\includegraphics[width=\textwidth]{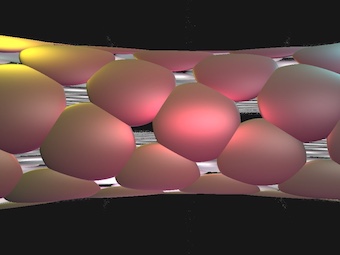}
   
\end{minipage}
~ 
\begin{minipage}[b]{0.32\textwidth} 
	\includegraphics[width=\textwidth]{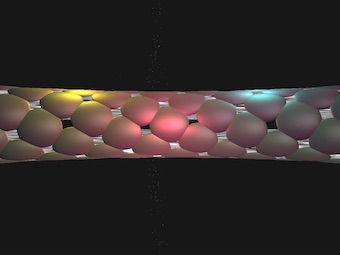}
\end{minipage}
\caption{Walking away from the $xy$-plane in $Sol$ along the $z$-axis at varying distances \href{https://vimeo.com/386630294}{(Video Link)}.}
\label{fig:dragon plane}
\end{figure}

Imagine that an observer stands in a hot-air balloon looking straight downwards (along the $-z$-axis) at the ground (the $xy$-plane).
As the balloon rises up, the plane appears as though it is ``rolled into a tube'' instead of extending infinitely far in all directions.
See \autoref{fig:dragon plane} (in the simulation, the plane has been tiled by small circles to give a sense of scale).

\begin{figure}
\centering
\includegraphics[trim=0 0 0 0,clip, width=0.8\textwidth]{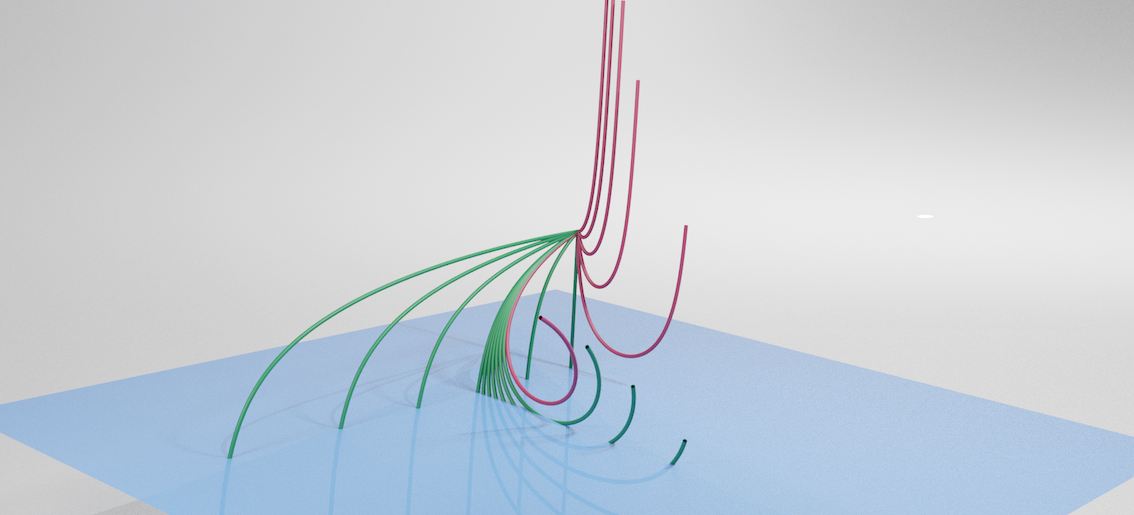}
\caption{A cluster of geodesics whose tangent vectors at time $0$ have negative dot product with $(0,0,1)$. Green rays hit the plane, red rays do not.}
\label{fig:dragon plane explained}
\end{figure}

\begin{figure}[h!tbp]
\centering
\begin{minipage}[b]{0.32\textwidth} 
	\includegraphics[width=\textwidth]{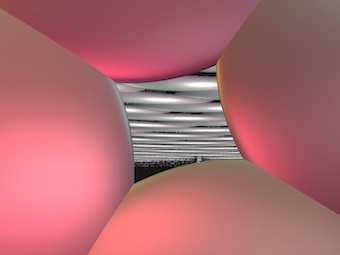}%
        	\subcaption{} 
\end{minipage}
~ 
\begin{minipage}[b]{0.32\textwidth} 
	\includegraphics[width=\textwidth]{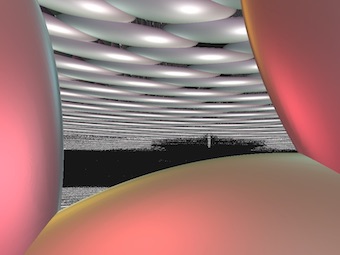}%
    \subcaption{} 
	\label{fig:through dragon plane - b}%
\end{minipage}
~ 
\begin{minipage}[b]{0.32\textwidth} 
	\includegraphics[width=\textwidth]{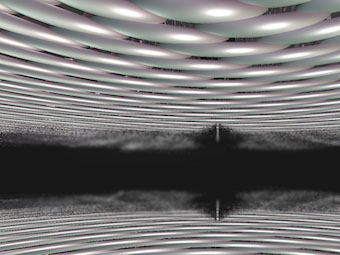}%
    \subcaption{} 
	\label{fig:through dragon plane - c}%
\end{minipage}
\caption{Walking forward through a plane in Sol.   (a) About to pass through the plane. 
(b) Passing through the plane, the backside is visible through the hole.  (c) After passing through, the backside of the plane still appears "in front". \href{https://vimeo.com/388598239}{(Video Link)}}
\label{fig:through dragon plane}
\end{figure}

Recall that as we are tracing our lines of sight along $\Sol$ geodesics, black points in the image correspond to directions for which the associated geodesics \emph{never} reach the $xy$-plane!
\autoref{fig:dragon plane explained} shows a cluster of geodesics starting above the $xy$-plane. 
Their tangent vectors all point downwards close to the $z$-axis.
Some of them, make a ``u-turn'' and head away in the other direction.
Those rays correspond to the dark area in the intrinsic view.
Another curious feature appears if the observer, still located above the $xy$-plane, looks through the holes of the circular tilling of the $xy$-plane.
It seems that there is something \emph{behind the plane}, see \autoref{fig:through dragon plane}.
Indeed, as shown on \autoref{fig:dragon plane explained}, some of the geodesic rays pointing downwards first hit the $xy$-plane, make a ``u-turn'', then hit the $xy$-plane again.
Thus, the light gray balls pictured in \autoref{fig:through dragon plane} correspond to the \emph{back side} of the $xy$-plane.
Imagine now that that the observer flies downwards through the $xy$-plane.
If their $z$-coordinate becomes negative, i.e. if the $xy$-plane is behind them, they will see the back side of the $xy$-plane \emph{in front of them}! 
See \autoref{fig:through dragon plane - c}.
This is caused again by the tendency of some geodesics to make a ``u-turn''.

\begin{figure}[h!tbp]
\centering
\begin{minipage}[b]{0.32\textwidth} 
	\includegraphics[width=\textwidth]{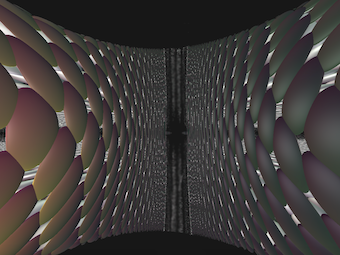}%
        	\subcaption{} 

\end{minipage}
~ 
\begin{minipage}[b]{0.32\textwidth} 
	\includegraphics[width=\textwidth]{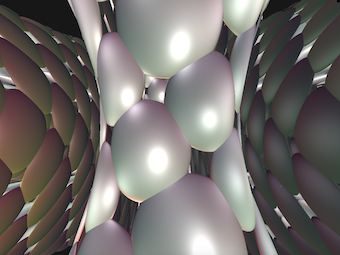}%
        	\subcaption{} 
\end{minipage}
~ 
\begin{minipage}[b]{0.32\textwidth} 
	\includegraphics[width=\textwidth]{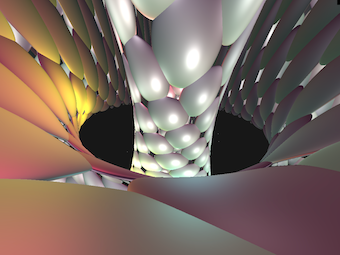}%
        	\subcaption{} 

\end{minipage}
\caption{Confusing views when sandwiched between two planes in Sol.  (a) Plane behind your head is visible in front of you, per above. (b) Drawing both the plane $z=-1$ behind you and the plane $z=1$ ahead of you. (c) A new perspective on the two planes. \href{https://vimeo.com/388610245}{(Video Link)}}
\label{fig: sandwich planes}
\end{figure}

For more confusing pictures, imagine that the observer stands at the origin $o = (0,0,0)$ ``sandwiched'' between two planes $\{z=-1\}$ and $\{z=1\}$.
The foreground plane appears as a rolled up tube, while the plane behind takes up a large fraction of the observer's forward-looking field of view, as per above.
Moving around between these planes offers a variety of interesting perspectives, see \autoref{fig: sandwich planes}.
In particular, when oriented the correct way the foreground plane may obscure the `vanishing line' of the background plane, giving the illusion that it is actually toroidal in shape.

\section*{Inside Compact Sol Manifolds}

The Lie group $X$ admits several \emph{uniform lattices}. A uniform lattice is a discrete subgroup $\Gamma \subset X$ whose corresponding quotient $M = \Gamma \setminus X$ is compact. 
An example of such a lattice is the subgroup $\Gamma$ generated by 
\begin{equation*}
    \gamma_1 =\left(\phi, -1,0\right) \quad
	\gamma_2 = \left(1, \phi,0\right) 
	\quad \text{and} \quad
	\gamma_3 = \left(0, 0,2 \ln \phi\right) ,
	\quad \text{where} \quad \phi = \frac {1 + \sqrt 5}2 \ \text{is the golden ratio}.
\end{equation*}

\begin{wrapfigure}{r}{0.42\textwidth}
	\centering
	\includegraphics[width=0.42\textwidth]{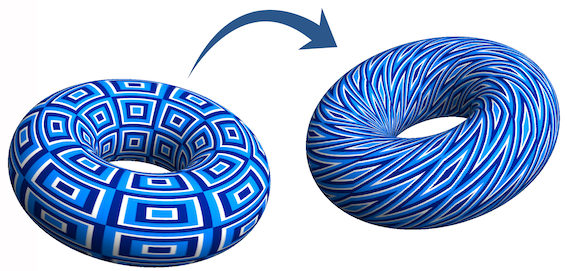}
\caption{The effect of an Anosov map on the torus.}
\label{fig:pA}
\end{wrapfigure}
The corresponding manifold $M$ is an example of one of the building blocks arising in Thurston's geometrization.
It has another interpretation, as follows.
Consider the matrix
\begin{equation}
	A = \left( 
	\begin{array}{cc}
		2 & 1 \\
		1 & 1
	\end{array}
	\right)
	\label{eqn: anosov}
\end{equation}
The action of $A$ on $\R^2$ preserves the integer points, i.e. $\Z^2 \subset \R^2$.
Hence $A$ induces a homemorphism $f$ of the two-dimensional torus $T$ (more precisely an \emph{Anosov homemorphism}) see \autoref{fig:pA}.
The \emph{mapping torus of  $T$ by $f$},  denoted by $T_f$, is the quotient of $T \times [0,1]$ by the equivalence relation which identifies each point $(x,1)$ with $(f(x),0)$.
We claim that $T_f$ can be endowed with a riemannian metric so that $M$ and $T_f$ are isometric.

\begin{figure}[h!tbp]
\centering
\begin{minipage}[b]{0.32\textwidth} 
	\includegraphics[width=\textwidth]{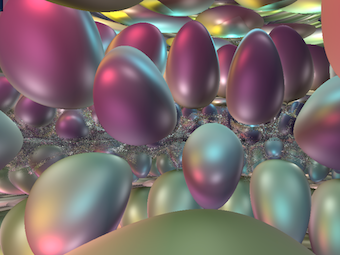}%
        	\subcaption{} 

\end{minipage}
~ 
\begin{minipage}[b]{0.32\textwidth} 
	\includegraphics[width=\textwidth]{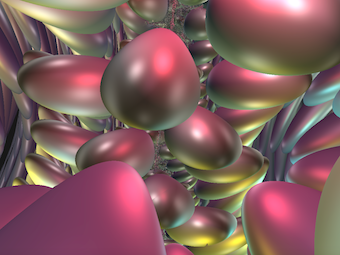}%
        	\subcaption{} 
\end{minipage}
~ 
\begin{minipage}[b]{0.32\textwidth} 
	\includegraphics[width=\textwidth]{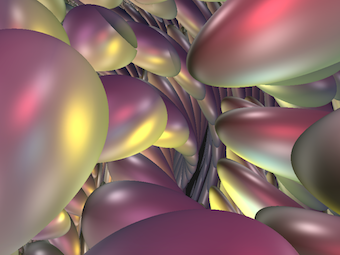}%
        	\subcaption{} 

\end{minipage}
\caption{A lattice in Sol.  (a) Looking in the $z$ direction (mapping torus direction).  (b) Looking in direction of $x$ axis in Sol, and (c) the direction $x=y$. \href{https://vimeo.com/388610289}{(Video Link)} }
\label{fig:lattice balls}
\end{figure}

Let us give a glimpse of this identification.
The action of $A$ on $\R^2$ stretches/compresses the plane in the direction of its eigenvectors.
Identifying these with the $x$- and $y$-axes in Sol shows the tori $T\times \{c\} \subset T_f$ to be the quotients of horizontal planes by the translations below
\begin{equation}
	g_1 \colon (x,y)\mapsto(x+\phi,y-1), 
	\quad
	g_2 \colon (x,y)\mapsto(x+1,y+\phi).
	\label{eqn:translation torus}
\end{equation}
These translations correspond to the elements $\gamma_1$ and $\gamma_2$ given above.
By construction, the matrix $A$ is diagonal in the $(x,y)$ coordinates. 
Hence the map $f$ sends $(x,y)$ to $(\phi^2 x, \phi^{-2}y)$ which corresponds to the action of $\gamma_3$.

We may use our technique - rendering images by tracing light rays outward from each pixel along geodesics - to visualize the interior of $M$.
Nevertheless, since $M$ is compact, any light ray ``wraps around'' the space many times. 
Consequently an inhabitant of $M$ could see infinitely many copies of the same object.
\autoref{fig:lattice balls} shows several inner views of $M$, which we populated with a \emph{single} ball.
Equivalently, it can be understood as a view of Sol where we positioned a ball at every lattice point $\gamma o$, for $\gamma \in \Gamma$.

\begin{figure}[h!tbp]
\centering
\begin{minipage}[b]{0.3\textwidth} 
	\includegraphics[width=\textwidth]{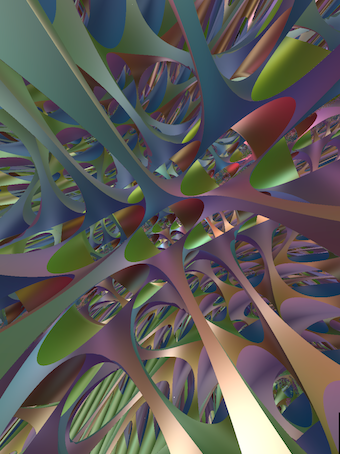}%
        	\subcaption{} 
        	\label{fig:2a}
\end{minipage}
~ 
\begin{minipage}[b]{0.3\textwidth} 
	\includegraphics[width=\textwidth]{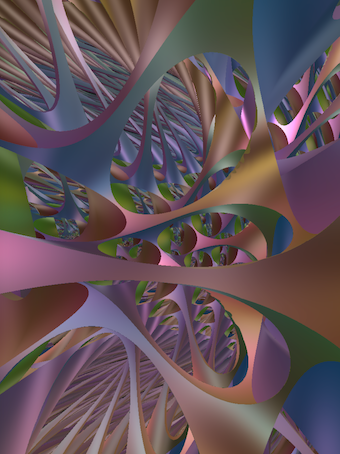}%
        	\subcaption{} 
        	\label{fig:2b}
\end{minipage}
~ 
\begin{minipage}[b]{0.3\textwidth} 
	\includegraphics[width=\textwidth]{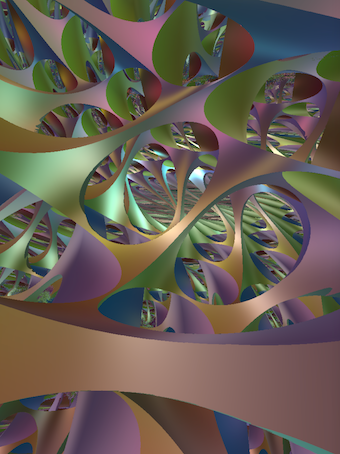}%
        	\subcaption{} 
        	\label{fig:2c}
\end{minipage}
\caption{Inside the Sol manifold with fundamental group drawn above.  (a) Looking in the $z$ direction.  (b) Looking in direction of $x$ axis in Sol, and (c) the direction $x=y$.
\href{https://vimeo.com/386628225}{(Video Link)}}
\label{fig:lattice pillar}
\end{figure}

\autoref{fig:lattice pillar} shows another point of view of the lattice $\Gamma$.
We join by a ``pillar/beam'' any two points $x$, $x'$ in $\Gamma o \subset X$ such that $x' = \gamma_i x$ for some $i \in \{1,2,3\}$.
These pillars stake out a fundamental domain of the action of $\Gamma$ on $X$ and its translates.
Alternatively the picture can be interpreted as follows.
Let $x_o$ be the image in $M$ of the origin $o$ of $X$.
The lattice $\Gamma$ is also the fundamental group of the quotient $M$.
In particular each $\gamma_i$ is represented by a simple closed curve $c_i$ based at $x_o$.
\autoref{fig:lattice pillar} shows what an observer living in $M$ sees if we thicken the curves $c_i$ into tubes.

\section*{Implementation Details}

We built our virtual reality Sol-simulator by adapting the technique of \emph{ray-marching} to non-euclidean homogeneous spaces.
This differs from previous work of the authors \cite{NEVR1,NEVR2}  which compute intrinsic views by pulling objects back to the tangent space via partial inverses of the riemannian exponential map.
Such a method cannot easily be applied for Sol, as the exponential map is far from being one-to-one.
In contrast, ray-marching, much like ray-tracing, works by flowing outwards from the screen along geodesics into the scene, and upon intersecting an object computing the color of the relevant pixel using material properties of that object, and the location / direction of light sources.
The precise implementation details are the subject of a forthcoming paper, and the code (currently a work in progress) is available on GitHub \cite{github}.
Other simulations of Sol geometry include the work of Berger \cite{berger},  ZenoRogue \cite{ZenoR} (both inverse-exponential and ray-tracing implementations) and MagmaMcFry \cite{MagmaMF} (ray-marching).

\paragraph{Producing a fixed image.}
Imagine that the observer stays at a fixed position in Sol without moving or rotating.
In order to compute the image she would see, we require a few geometric ingredients: a means of computing the geodesic flow in $M$, and a \emph{signed distance function} measuring the distance in $M$ from each point to the nearest object in the scene (this tells us how far we can safely flow along a geodesic without hitting an object).
As mentioned above, the geodesic flow in Sol can be solved explicitly using Jacobi's elliptic/zeta functions \cite{Troyanov:1998aa}.
Note that since Sol has no continuous symmetries fixing the origin $o$, it is not possible to reduce the dimension of the problem.

Distance functions are extremely difficult to compute explicitly in Sol.
This comes from the fact that many geodesics do not globally minimize distance.
Moreover the form of the geodesic flow is complicated.
Nevertheless, the distance function to a plane $\{ z = c\}$ is rather easy: as the vertical geodesics $t \to (x,y,z+t)$ are minimizing, the signed distance from the point $(x,y,z)$ to the $xy$-plane is simply the $z$ coordinate.
This allows us to accurately render planes of this form in Sol.
Other objects are rendered using a ``fake distance function'' which approximates the true distance.

\paragraph{Moving in the Space.}
Imagine now that the observer wants to walk and explore the space around her.
In order to render the image she would see, we need a way to compute her position and facing.
We decided that a straight displacement (if the user holds down the forward key) should move the observer along the geodesic whose tangent vector is given by the direction she is looking, while her orientation around the geodesic is updated using parallel transport.
In contrast to isotropic spaces (such as the euclidean space) there is no $1$-parameter group of isometries realizing both displacement and the parallel transport.

In our software, we encode the displacement using elements of Sol (recall that Sol acts freely transitively on itself by isometries).
As for the parallel transport operator, we use an idea explored by Grayson \cite{Grayson:1983aa}.
Let $c \colon \R_+ \to X$ be a geodesic starting at the origin $o$.
The parallel transport from $c(0)$ to $c(t)$ along $c$ is an isometry $T(t) \colon T_{c(0)}X \to T_{c(t)}X$.
To avoid any confusion we denote by $L_{c(t)}$ the element of Sol (seen as an isometry of $X$) sending the origin to $c(t)$.
This can be used to pull back the parallel transport to the tangent space at the origin.
More precisely we focus on the operator of $Q(t) \colon  T_oX \to T_oX$ defined by $Q(t)  = dL_{c(t)}^{-1} T(t)$.
Identifying $T_oX$ with $\R^3$, the operator $Q(t)$ is a matrix in ${\rm SO}(3)$ satisfying
\begin{equation*}
	\dot Q + BQ = 0, 
	\quad \text{with} \quad
	B = \left(
		\begin{array}{ccc}
			0 & 0 & -u_x \\
			0 & 0 & u_y \\
			u_x & -u_y & 0
		\end{array}
	\right), 
\end{equation*}
where $u = (u_x,u_y,u_z)$ is defined by $u(t) = dL_{c(t)}^{-1}\dot c(t)$.
In practice, all computations here are made with the Runge-Kutta method.

\paragraph{Quotient Manifolds.}
To ray-march in the Sol manifold $M = \Gamma \setminus X$ depicted above, we identify $M$ with a fundamental domain $D$ for the action of $\Gamma$ on $X$.
We defined an algorithm, so that every time a light ray escapes $D$ it is ``teleported'' back into $D$ using an element of $\Gamma$.
In this way the orbits of our flow in $D$ map to geodesics in $M$.
In the extrinsic model $(\R^3, ds^2)$ of Sol, the fundamental domain $D$ has the form $D_0 \times [0,2\ln \phi)$, where $D_0$ is a fundamental domain for the action of the translations $g_1$ and $g_2$ on $\R^2$, see Equation (\ref{eqn:translation torus}).
The ``teleportation'' works as follows.
Assume that the point $p = (x,y,z)$ does not belong to $D$.
Translating $p$ by a suitable power of $\gamma_3$ we can make sure that $p$ belongs $\R^2 \times [0,2\ln \phi)$.
Then, in either order (as they commute) we iteratively apply $\gamma_1$ and $\gamma_2$ (which do not affect $z$) until the point has been brought back into $D$.


\section*{Summary and Future Work}

In summary, this project has produced a real-time, intrinsic and geometrically correct rendering engine for Sol geometry and its compact quotients, that has the ability to take movement input from either a keyboard or headset, and has the ability to render images from two viewpoints simultaneously (for stereoscopic vision).
However, this is still very much a work in progress.
The list of features still under development include
\begin{enumerate}
\item a good approximation of the $\Sol$ distance function to a point, for intrinsic rendering of geodesic spheres, 
\item tracking down and cleaning up the sources of noise in the numerical approximation to Jacobi functions required for the geodesic flow, 
\item computational speedup to allow real-time rendering in high-definition, and 
\item given two points $\{p,q\}$, a procedure for calculating the tangent vectors at $p$ which correspond to geodesics reaching $q$ (for accurate lighting considerations).
\end{enumerate}


\vspace{-0.25cm}
\section*{Acknowledgements}

This material is based upon work supported by the National Science Foundation under Grant No. DMS-1439786 while the authors were in residence at the ICERM in Providence, RI, during the semester program \emph{Illustrating Mathematics}.
We are thankful to Brian Day in getting things synched up with the virtual reality headsets.
We are additionally grateful to many others at ICERM for interesting conversations about Sol geometry, including Rich Schwartz, Matei Coiculescu and Jason Manning.
The first author acknowledges support from the Agence Nationale de la Recherche under Grant \emph{Dagger} ANR-16-CE40-0006-01 as well as the \emph{Centre Henri Lebesgue} ANR-11-LABX-0020-01. The second author is grateful to support from the National Science Foundation DMR-1847172. The third author was supported in part by National Science Foundation grant DMS-1708239.
This project is indebted to a long history of previous work.  
It is a direct descendant of the hyperbolic ray-marching program created by Nelson, Segerman, and Woodard~\cite{RSC}, which itself was inspired by previous
work in $\mathbb H^3$ and $\mathbb H^2\times \mathbb E$ by  Hart, Hawksley,  Matsumoto, and Segerman \cite{NEVR1,NEVR2}, all of which aim to expand upon the excellent work of Weeks in \emph{Curved Spaces}~\cite{Weeks}.


\vspace{-0.25cm}

    
{\setlength{\baselineskip}{13pt} 
\raggedright				

} 
   
\end{document}